\documentclass[11pt]{article}

\usepackage{amsfonts,amssymb}
\usepackage{graphicx}

\newcommand{\Z}{\mathbb{Z}}
\newcommand{\R}{\mathbb{R}}
\newcommand{\Q}{\mathbb{Q}}
\newcommand{\F}{\mathbb{F}}

\newcommand{\id}{\mbox{\rm id}}

\newcommand{\vcob}{{\bf VCob}}

\newcommand{\mcob}{{\bf Cob}_{F \times I}}
\newcommand{\vd}{\rm{deg}}

\def\id{{|}}

\newtheorem{theorem}{Theorem}[section]
\newtheorem{corollary}[theorem]{Corollary}
\newtheorem{lemma}[theorem]{Lemma}
\newtheorem{proposition}[theorem]{Proposition}

\newtheorem{definition}[theorem]{Definition}

\newtheorem{example}[theorem]{Example}

\newtheorem{remark}[theorem]{Remark}

\input{epsf.sty}

\begin{document}

\title{Frobenius Modules and  Essential Surface Cobordisms{\footnote{MRCN:16W30;57M99;55u99;18A99}}}

\author{
J. Scott Carter \\
University of South Alabama \\
\and
Masahico Saito\\
University of South Florida 
}

\maketitle

\begin{abstract}
An algebraic system is proposed that represent surface cobordisms 
in thickened surfaces. Module and comodule structures over Frobenius algebras are used for 
representing essential curves. The proposed structure gives a unified algebraic view 
of states of categorified Jones polynomials in thickened surfaces and virtual knots.
Constructions of such system are presented. 
\end{abstract}

\section{Introduction}

In  this article, we propose a formulation of an algebraic structure 
that describes 
surface cobordisms in thickened surfaces 
 that have 
both inessential  
and
essential circles. 
Thus the structure we propose  is a refinement 
of   a   $(1+1)$-TQFT.
Our motivation comes from 
 the differentials of  generalizations of Khovanov homology~\cite{Kh99} 
defined  in \cite{APS,TT} 
 for thickened surfaces and
those 
for virtual knots~\cite{IT,Man}. 
Although the use of a Frobenius algebra structure is explicit in \cite{IT,TT}, 
we propose 
a more detailed distinction between the vector spaces 
that are assigned to trivial circles  
and  
essential 
circles in  thickened surfaces. 
In this way, we will construct an algebraic system that provides a unified view 
of the states and differentials used in these theories.  
We propose that essential circles can be studied by means of a module and comodule structure over a Frobenius algebra.

The well-definedness of the differential ($d^2=0$) in the Khovanov homology for classical knots
depends mainly on a  $(1+1)-$TQFT structure 
and therefore upon a Frobenius algebra~\cite{Kh06}. 
Variations of a  $(1+1)-$TQFT structure have been studied (\cite{Carmen,Lauda}, for example)
for refinements of the Khovanov homology.  
Studies of generalizations of TQFT to surfaces in $3$-manifolds were suggested also 
by C. Frohman. 
In \cite{APS}, 
the authors adapt the  description of \cite{Viro} and define differentials in a combinatorial fashion
by using signs ($\pm)$ and enhanced states to define their differentials.
 In particular, they don't use an explicit Frobenius algebra. 
On the other hand, an algebraic formulation creates the advantage of enabling systematic generalizations and streamlining the proofs of well-definedness~\cite{Kh06}. 
Meanwhile, in \cite{IT,TT} a Frobenius algebra is used to generalize Khovanov homology to 
virtual knots 
and thickened surfaces, respectively. 
Herein, we provide a single algebraic approach that  envelopes both theories, and we provide constructions and examples for our approach.

The paper is organized as follows.
In Section~\ref{prelimsec}, 
we review necessary materials  and establish notation.   
In Section~\ref{defsec}, we define the algebraic structure 
called  ``commutative Frobenius pairs" and present examples.
We show in Section~\ref{TQFTsec}  that these structures naturally arises in 
$(1+1)$-TQFTs. 
In Section~\ref{constructsec}, 
we provide 
new  
methods of constructions.

\section{Preliminaries} \label{prelimsec}

\subsection{Frobenius algebras and their diagrams} 
Frobenius algebras are assumed to be as described  in \cite{Kock}, Section 1.3, 
and we give a brief summary here.  
A {\it Frobenius algebra} is an (associative) algebra (with
multiplication $\mu: A \otimes A \rightarrow A$ and unit $\eta: k
\rightarrow A$) over a unital commutative ring 
$k$ with a nondegenerate associative
pairing $\beta:  A \otimes A \rightarrow k$. 
The  pairing $\beta$ is also expressed by
$\langle x|y\rangle=\beta(x \otimes y)$ for $x, y \in A$, and it
is {\it associative}
in the sense that
$\langle xy | z\rangle=\langle x|y
z\rangle$ for any $x,y,z \in A$.

\begin{figure}[htb]
\begin{center}
\includegraphics[width=4.5in]{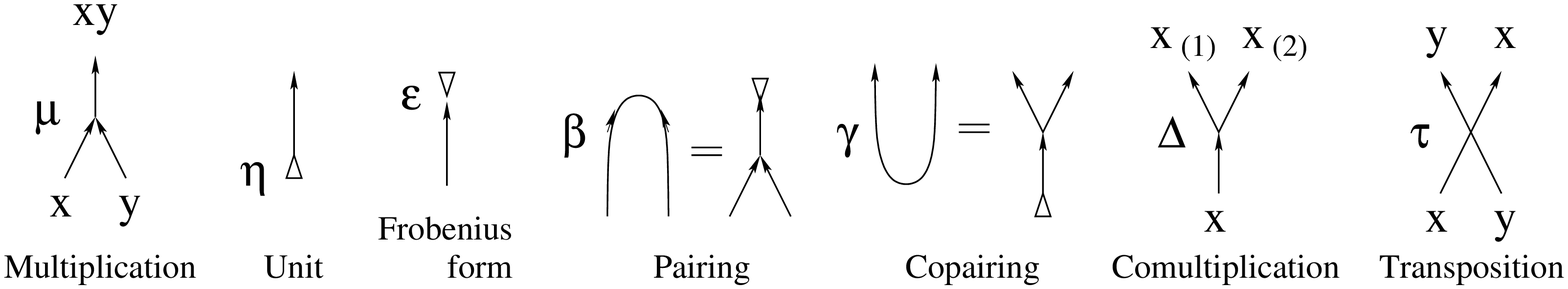}
\end{center}
\caption{Diagrams for Frobenius algebra maps}
\label{frodiag}
\end{figure}

\begin{figure}[htb]
\begin{center}
\includegraphics[width=4in]{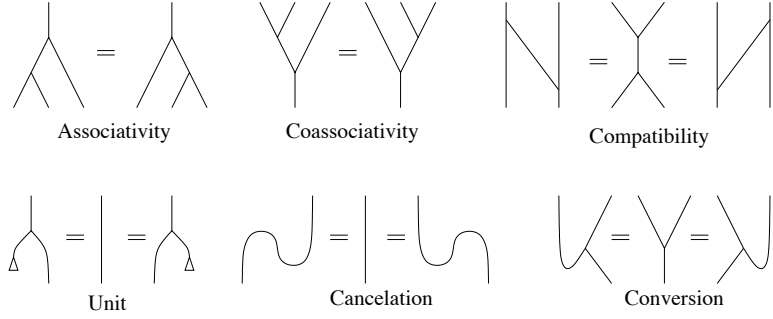}
\end{center}
\caption{Equalities among Frobenius algebra maps }
\label{froaxioms}
\end{figure}

A Frobenius algebra $A$ has a linear functional $\epsilon: A
\rightarrow k$,
called the {\it Frobenius form},
or a {\it counit}, 
 such that the kernel
 contains no nontrivial  left ideal.
 It is defined from $\beta$ by $\epsilon(x)=\beta(x \otimes 1)$,
and conversely, a Frobenius form gives rise to a nondegenerate
associative pairing $\beta$ by $\beta(x \otimes y)=\epsilon(xy)$,
for $x, y \in A$.
A Frobenius form has a unique copairing
$\gamma: k \rightarrow A \otimes A$
characterized by
$$ (\beta \otimes |)(| \otimes \gamma) = |
 = (| \otimes \beta)(\gamma \otimes |) , $$
 which we call the
 {\it cancelation} of $\beta $ and $\gamma$.
  See the middle entry in the bottom row of Fig.~\ref{froaxioms}.
 Here and below, we denote by $|$ the identity homomorphism on the algebra. This notation will distinguish this function from the identity element $1=1_A=\eta(1_k)$ 
 of the algebra that is the image of the identity of the ground 
 ring. 
A Frobenius algebra $A$ determines a coalgebra structure with
$A$-linear (coassociative) comultiplication and the counit
defined using the Frobenius form. The comultiplication  $\Delta: A
\rightarrow A \otimes A$ is defined by
$$ \Delta  =  
 (\mu \otimes |)(| \otimes \gamma) =
 (| \otimes \mu)(\gamma \otimes |). $$
The multiplication and comultiplication
satisfy the following equality:
$$ \Delta \mu = (\mu \otimes |)(| \otimes \Delta )=(| \otimes \mu)(\Delta \otimes |)$$
which we call the {\it Frobenius compatibility condition}.
In Fig.~\ref{frodiag}, diagrammatic conventions of 
various maps that appear for Frobenius algebras are depicted. 
The diagrams are read from bottom to top, and each line segment represents a tensor factor of $A$.
In Fig.~\ref{froaxioms}, the axioms and relations among compositions 
of these maps are represented by these diagrams.

A Frobenius algebra is {\it commutative} if it is commutative as
an algebra.   It is known (\cite{Kock} Prop. 2.3.29) that
a Frobenius algebra is commutative if and only if it is
cocommutative as a coalgebra.
The map $\mu \Delta $ of a Frobenius algebra
 is called the
 {\it handle operator,} and
 corresponds to multiplication by a central element called the {\it handle element}
$\delta_h =\mu \Delta (1)$ 
 (\cite{Kock}, page 128).

\begin{example} \label{univexample}
{\rm
This example appears in universal Khovanov homology \cite{Kh99}.
Let $A=\Z[X, h, t]/(X^2 -hX -t) $, with unit $\eta(1)=1$, counit $\epsilon(1)=0$, $\epsilon(X)=1$,
multiplication being the polynomial multiplication, 
comultiplication
defined by 
\begin{eqnarray*}
\Delta(1)&=& 1 \otimes X + X \otimes 1 - h \ 1 \otimes 1\\
\Delta(X)&=& X \otimes X +t \ 1 \otimes 1.  
\end{eqnarray*}

}\end{example}

\subsection{Modules and comodules}

In this paper we focus on commutative, cocommutative algebras and bimodules
and bicomodules, and assume the following conditions.
Let $k$ be a unital commutative ring. 
A commutative algebra is a ring $A=(A, \mu, \eta)$, 
 that is 
a $k$-module $A$ 
 $k$-linear 
multiplication $\mu: A \otimes A \rightarrow A$
is associative, with 
its $k$-linear unit map  denoted by $\eta=\eta_A: k \rightarrow A$.
We also use the notations $1=1_k$, $1_A=\eta (1)$.  
By a {\it commutative bimodule} $E$ over $A$, we mean that 
$E$ is an $A$-bimodule, and 
the left and right actions
$\mu= \mu_{A, E}^E: A \otimes E \rightarrow E$ 
(denoted by $a \otimes x \mapsto ax=\mu(a \otimes x$) ), 
$\mu= \mu_{E, A} ^E : E \otimes A \rightarrow E$ 
(denoted by $x \otimes a \mapsto xa=\mu(x \otimes a$) )  
satisfy $ax=xa$  for any $a \in A$, $x \in E$.
Recall that the module conditions include  that 
$a(bx)=(ab)x$ and $1_A x = x$ for 
all $a, b \in A$ and $x \in E$.

\begin{figure}[htb]
\begin{center}
\includegraphics[width=4in]{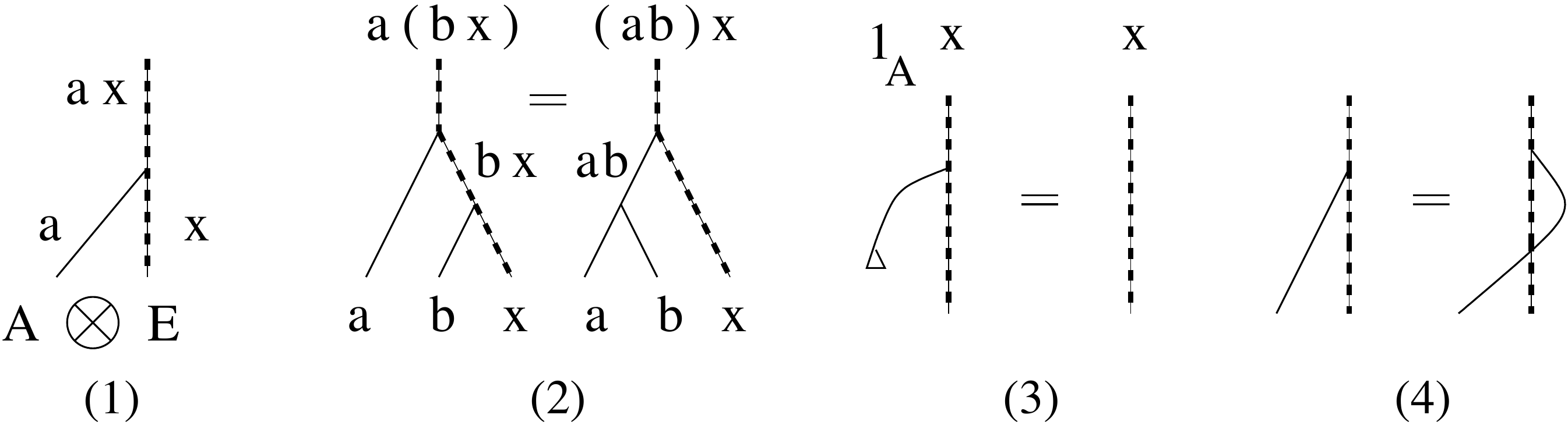}
\end{center}
\caption{Bimodule  maps and their relations}
\label{module}
\end{figure}

In diagrams, we represent the $A$-module  $E$ by thick dotted lines as depicted in Fig.~\ref{module}.
In the figure, the following maps and formulas are depicted for $a, b \in A$ and $x \in E$:
(1) the action $a \otimes x \mapsto ax$, (2) $a(bx)=(ab)x$, (3) $(1_A) \cdot x=x$, and
(4) $ax=xa$. 

Let $V$, $W$ be free $k$-modules of finite rank for a unital commutative ring $k$.
Denote by  
$\tau_{V,W}: V \otimes W \rightarrow W \otimes V$ for $k$-modules $V, W$ the 
$k$-linear map induced from 
the transposition
 $\tau(\sum v \otimes w)= \sum w \otimes v $.

Let $B=(B, \Delta, \epsilon)$ be a cocommutative coalgebra over $k$
with coassociative comultiplication $\Delta$ and counit map $\epsilon: B \rightarrow k$.
By a {\it cocommutative bicomodule} $E$ over $B$, we mean that
$E$ is a bicomodule over $A$, and 
the coactions  ($k$-linear maps)
$\Delta=\Delta_E^{B, E} : E \rightarrow B \otimes E$
 (denoted by $ x \mapsto \sum x_{(0)} \otimes x_{(1)}$), 
 $\Delta=\Delta_E^{E,B} : E \rightarrow E \otimes B$
  (denoted by $ x \mapsto \sum x_{(1)}'  \otimes x_{(0)}'$), 
satisfy
$\tau (\sum x_{(0)} \otimes x_{(1)}) =\sum x_{(1)}'  \otimes x_{(0)}'$.
for any $a \in A$, $x \in E$.
The corresponding figures are upside-down diagrams of Fig.~\ref{module}.

\section{Frobenius pairs} \label{defsec}

\begin{figure}[htb]
\begin{center}
\includegraphics[width=3in]{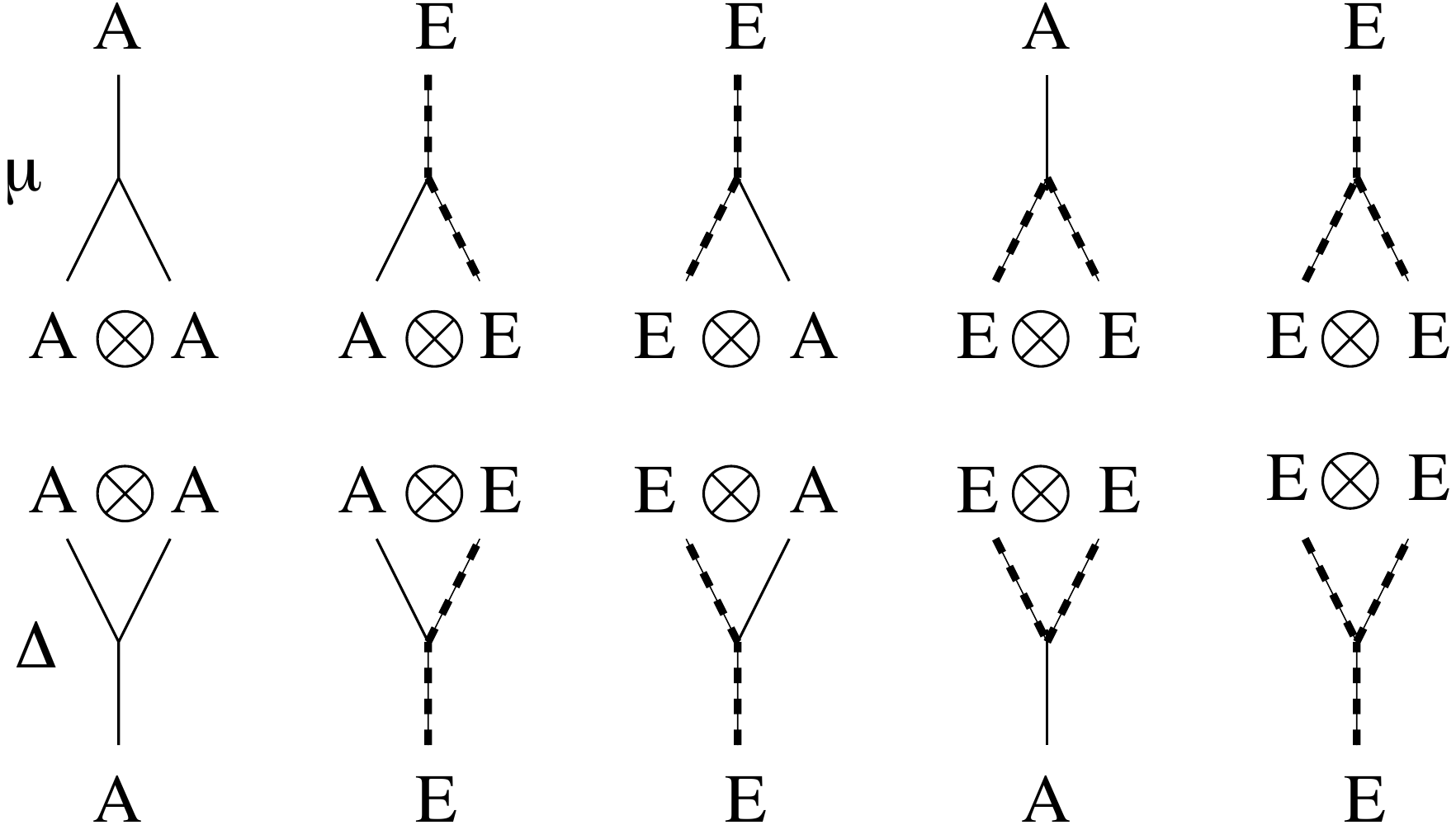}
\end{center}
\caption{Generating maps}
\label{gens}
\end{figure}

\begin{figure}[htb]
\begin{center}
\includegraphics[width=2.5in]{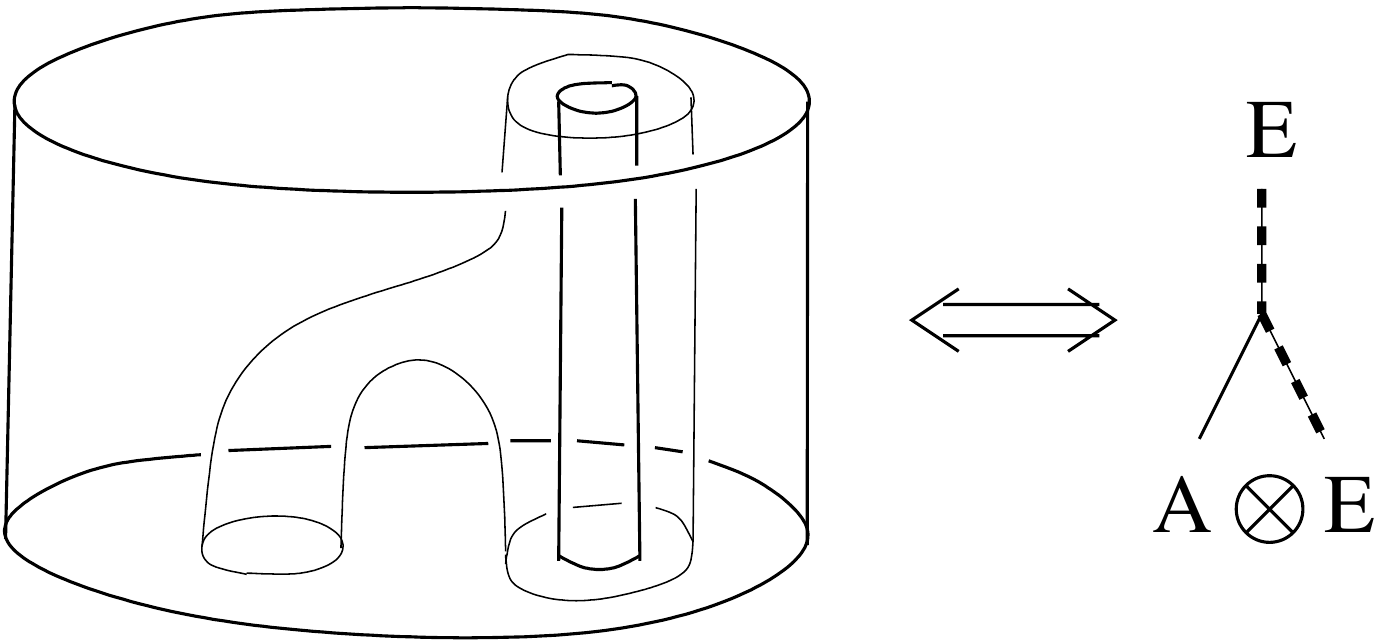}
\end{center}
\caption{A cobordism with essential curves and module action}
\label{essentialcob}
\end{figure}

In this section, we give 
the 
definitions of 
the 
algebraic structures that 
are studied in this paper. 
Our diagrammatic convention to represent these maps
 in the definition below 
is depicted in Fig.\ref{gens}. 
The ring $A$ and 
the  $A$-module  
$E$ are represented by solid and dotted lines, 
respectively, and variety of multiplication and comultiplication 
are depicted by trivalent vertices, read from bottom to top. 
We note that the possibilities of trivalent vertices are 
summarized by saying that the dotted line does not end at a trivalent vertex, while 
a solid line can. 
The definition below is motivated from  surface cobordisms, 
and the correspondence is exemplified in Fig.~\ref{essentialcob}.
Briefly, the elements of $E$ are associated to essential curves in the surface cobordism and the elements of $A$ are associated to compressible curves.  
More details on the correspondence 
are given in Section~\ref{TQFTsec}.

\begin{definition} {\rm
A {\it commutative Frobenius  pair} $(A, E)$ is 
defined as follows.

\smallskip

\noindent
(i) 
$A=( \mu_A, \Delta_A, \eta_A, \epsilon_A)$ 
is a commutative Frobenius algebra over $k$
with multiplication $ \mu_A$, comultiplication $\Delta_A$,
unit $\eta_A$ and counit $\epsilon_A$.

\smallskip

\noindent
(ii) $E$ is an $A$-bimodule and $A$-bicomodule, 
with the same  right and left 
  actions and coactions.

\smallskip

\begin{figure}[htb]
\begin{center}
\includegraphics[width=2in]{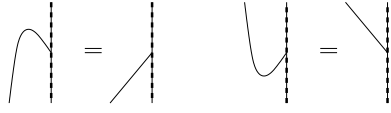}
\end{center}
\caption{Cancelation relations}
\label{cancel}
\end{figure}

\noindent
(iii) The action and coaction satisfy the canceling conditions with the pairing and copairing
as follows:
$  (\beta \otimes \id) (\id \otimes \Delta_E^{A,E})= \mu_{A,E}$,
$ (\id \otimes  \mu_{A,E}^E)(\gamma \otimes \id)=\Delta_E^{A,E}$. 
The situation is depicted in Fig.~\ref{cancel}.

\noindent
(iv)
$E$ has  
an associative, commutative  multiplication $ \mu_E: E \otimes E \rightarrow E$
and  a coassociative commutative comultiplication 
 $\Delta_E:  E \rightarrow  E \otimes E$,
 that are $A$-bimodule and $A$-bicomodule maps, 
 such that 
the maps $ \mu_E$ and $\Delta_E$ satisfy the compatibility condition:
$$(\id \otimes  \mu_E)(\Delta_E \otimes \id )
= \Delta_E  \mu_E
= ( \mu_E \otimes \id )(\id \otimes \Delta_E) .  $$
The diagrams for these relations are the same as associativity for $\mu_A$, 
except that all segments are dotted. 

\smallskip

\begin{figure}[htb]
\begin{center}
\includegraphics[width=3in]{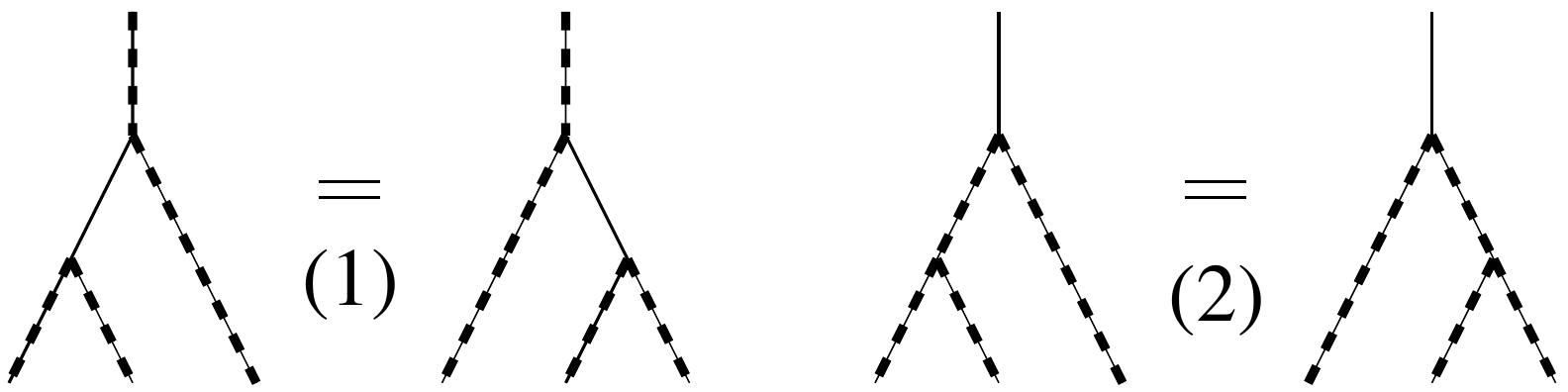}
\end{center}
\caption{Analogs of associativity}
\label{EEA}
\end{figure}

\noindent 
(v) There are $A$-bimodule, $A$-bicomodule maps
$ \mu_{E,E}^{A}
:  E \otimes E \rightarrow A$ and 
$ \Delta_{A}^{E,E} 
: A \rightarrow E \otimes E$
that are associative and coassocative, respectively:
\begin{eqnarray*}
 \mu_{A,E}^E (  \mu_{E,E}^{A}  \otimes \id_E) =  \mu_{E,A} ^E(  \id_E  \otimes  \mu_{E,E}^{A} ) &:& 
 E \otimes  E \otimes E \rightarrow E , \label{v1} \\
  \mu_{E,E}^{A} (  \mu_{E}  \otimes \id_E) =  \mu_{E,E}^{A} (  \id_E  \otimes  \mu_{E} ) &:& 
 E \otimes  E \otimes E \rightarrow A ,   \label{v2}  \\
(\Delta_{A}^{E,E} \otimes \id_E) \Delta_E^{A,E} = (\id_E \otimes \Delta_{A}^{E,E}) \Delta_E^{E,A}
&:&  E \rightarrow E \otimes E \otimes E , \label{v3} \\
(\Delta_E \otimes \id_E) \Delta_{A}^{E,E}= (\id_E \otimes \Delta_E) \Delta _{A}^{E,E}
&:&  A \rightarrow E \otimes E \otimes E . \label{v4} 
\end{eqnarray*}
The first two are depicted in Fig.~\ref{EEA}, and the last two are their upside-down diagrams.

\smallskip

\begin{figure}[htb]
\begin{center}
\includegraphics[width=3in]{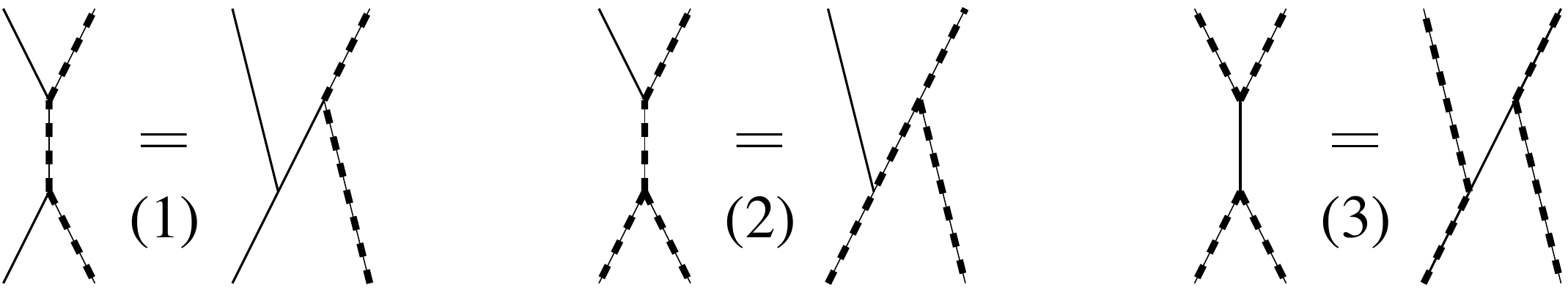}
\end{center}
\caption{Compatibility conditions} 
\label{compati}
\end{figure}

\noindent 
(vi) These maps satisfy the compatibility condition depicted in Fig.~\ref{compati}, 
that are analogs of the compatibility condition of  multiplication and comultiplication of Frobenius algebras.
Specifically,
\begin{eqnarray*}
\Delta_E^{A,E}   \mu_{A,E}^E &=& (  \id_A \otimes  \mu_{A,E}^E )   ( \Delta_A  \otimes  \id_E ) , \\
\Delta_E^{A,E}   \mu_{E,E}^E &=& 
(  \id_A \otimes  \mu_{E,E}^E )   ( \Delta_E^{A,E}  \otimes  \id_E ) 
\\
\Delta_A^{E,E}   \mu_{E,E}^A &=& 
(  \id_E \otimes  \mu_{A,E}^E )   ( \Delta_E^{E,A}  \otimes  \id_E )  .
\end{eqnarray*}
We also include the equalities represented by upside-down and mirror image 
diagrams of  Fig.~\ref{compati}.

\smallskip

\begin{figure}[htb]
\begin{center}
\includegraphics[width=4in]{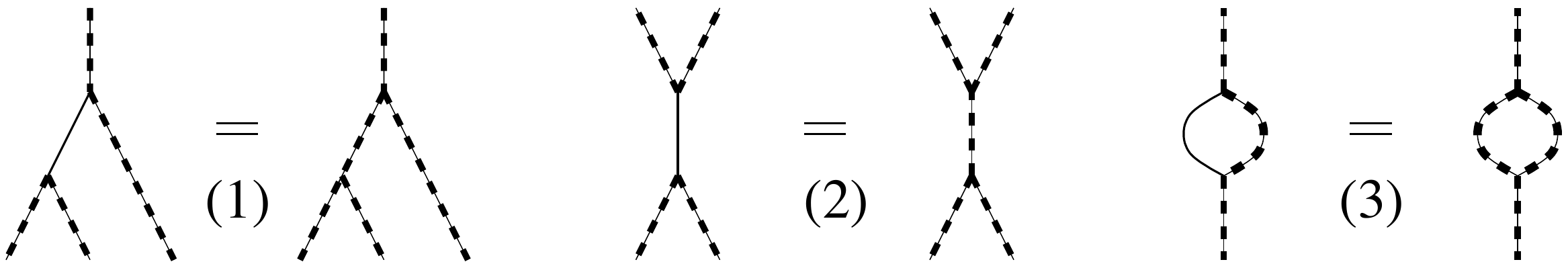}
\end{center}
\caption{Consistency conditions}
\label{samerel}
\end{figure}

\noindent 
(vii) The following relations, called {\it consistency conditions}, are satisfied: 
\begin{eqnarray*}
   \mu_{A,E}^E  ( \mu_{E,E}^A\otimes \id_E)  &=&    \mu_{E} (   \mu_{E} \otimes \id_E) , \\
  \Delta_A^{E,E}   \mu_{E,E}^A &=&\Delta_E   \mu_E ,\\
  \mu_{A,E}^E \Delta_E^{A,E} &=&   \mu_E \Delta_E. 
\end{eqnarray*}
These relations are depicted in  Fig.~\ref{samerel} (1), (2), and (3), respectively. 

} \end{definition}

\begin{figure}[htb]
\begin{center}
\includegraphics[width=5in]{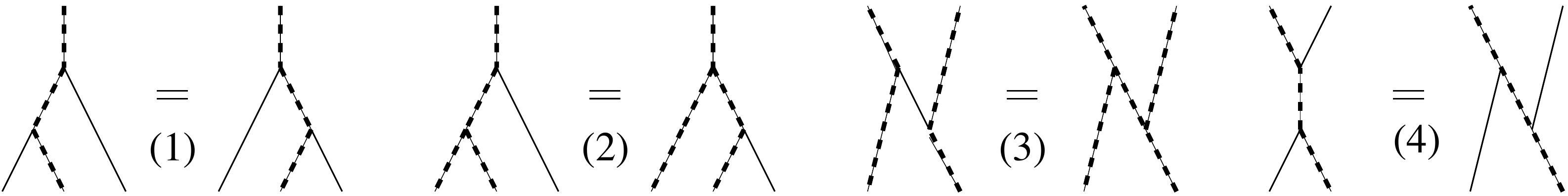}
\end{center}
\caption{Equalities that follow from definitions}
\label{maplemmas}
\end{figure}

It is left as an exercise to prove that the relations 
depicted in Fig.~\ref{maplemmas},
as well as their upside-down and mirror image diagrams, 
follow from the definitions.

\begin{figure}[htb]
\begin{center}
\includegraphics[width=5in]{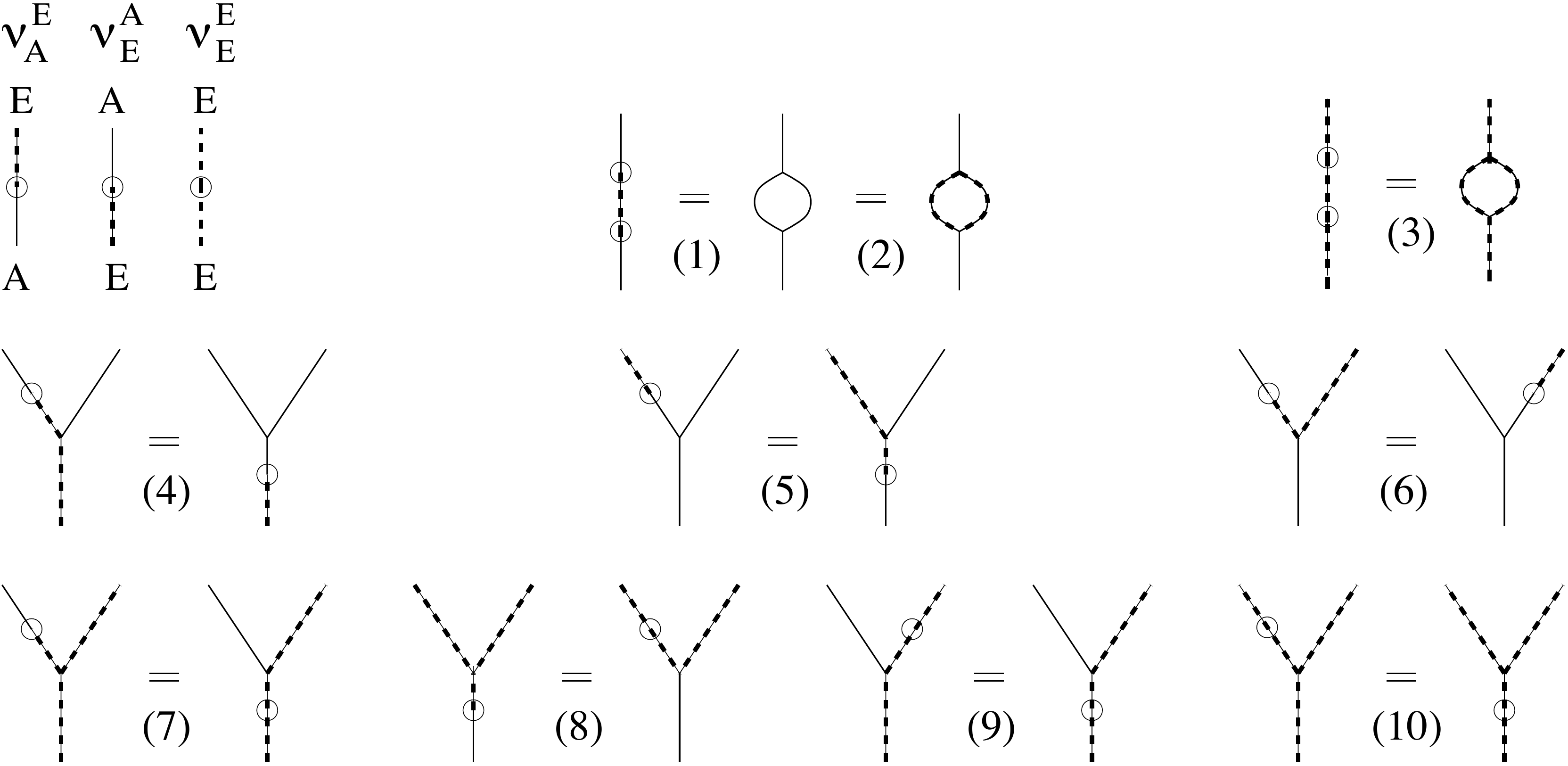}
\end{center}
\caption{M\"{o}bius maps and their relations}
\label{mobius}
\end{figure}

\begin{definition} {\rm 
Let $(A, E)$ be a commutative Frobenius pair. 
Three $k$-linear maps $\nu_{A}^{E}: A  \rightarrow   E$,
 $\nu_{E}^{A}: E  \rightarrow   A$ and 
  $\nu_{E}^{E}: E  \rightarrow   A$
are called  { \em M\"{o}bius maps } if they satisfy
the following conditions.
\begin{eqnarray*} 
\nu^A_E \nu^E_A 
= \mu_A \Delta_A = \mu_{E,E}^A \Delta_A^{E,E} :  A   \rightarrow A, & \ \ \  &   
  (\nu_{E}^{E} )^2 =  \mu_E \Delta_E : E \rightarrow E,  \\ 
  \mu_{E,A}^E  (\nu_{A}^{E} \otimes \id_E) = \nu_{A}^{E}  \mu_{A}
: A  \otimes E \rightarrow E, & \ \ \  &
(\nu^A_E \otimes \id_A)\Delta^{E,A}_E =\Delta_A\nu^A_E:E \rightarrow A\otimes A, 
 \\ 
   \mu_{A}  (\nu_{E}^{A} \otimes \id_A) = \nu_{E}^{A}  \mu_{A, E}^E
: E  \otimes A \rightarrow E, &  \ \ \  & (\nu^{E}_A \otimes \id_A)\Delta_A = \Delta^{E,A}_E\nu^E_A: A \rightarrow E \otimes A, \\ 
 \mu_E  (\nu_{A}^{E} \otimes \id_E) = \nu_{E}^{E}  \mu_{A,E}^E
: A  \otimes E \rightarrow E, &  \ \ \  &
(\nu^A_E \otimes \id_E) \Delta^{E,E}_A = (\id_A \otimes \nu^E_A)\Delta_A: A \rightarrow A\otimes E,
 \\ 
 \mu_E  (\nu_{A}^{E} \otimes \id_E) = \nu_{E}^{E}  \mu_{A,E}^E
: A  \otimes E \rightarrow E, &  \ \ \  & (\nu^A_E \otimes \id_E) \Delta_E =\Delta^{A,E}_E \nu^E_E:E \rightarrow A\otimes E, \\ 
\nu_{A}^{E}  \mu_E =\mu_{E,E}^A
: E  \otimes E \rightarrow A, &  \ \ \  &  \Delta_E \nu^E_E  =(\nu^E_E \otimes \id_E)\Delta^{E,E}_A:A \rightarrow E\otimes E,
\\ 
\mu_{A,E}^A (\id_A \otimes \nu_E^E )= \nu_{E}^{E}  \mu_{A,E}^E 
: A  \otimes E \rightarrow A, &  \ \ \  &  (\id_A \otimes \nu^E_E)\Delta^{A,E}_E  =\Delta^{A,E}\nu^E_E:E \rightarrow A\otimes E,
\\  
\mu_E (\nu_E^E \otimes \id_E )= \nu_E^E \mu_E
: E  \otimes E \rightarrow E, &  \ \ \  &  (\nu^E_A \otimes \id_E)\Delta_E  =\Delta_E \nu^E_E :A \rightarrow A\otimes E,
\\   
 \end{eqnarray*}
 Diagrams for the M\"{o}bius maps  $\nu_{A}^{E}$, $\nu_{E}^{A}$ and $\nu_{E}^{A}$
 are depicted in Fig.~\ref{mobius} at the top left.
 The required equalities are depicted in the same figure: items (1) through (3) are written at the top line of this equation array. The right-hand relations are depicted in (4) through (10) while the left-hand relations are the upside-down versions of these diagrams. 
} \end{definition}

We abuse the notation for $m$ and $\Delta$ by dropping subscripts
if no confusion occurs. 
Thus, $m : (A \otimes A)  \oplus (A \otimes E )\oplus (E \otimes E)
\rightarrow A \oplus E$ represents 
$ \mu_A +   \mu_{A,E}^E +  \mu_E$ and 
$\Delta : A \oplus E \rightarrow 
(A \otimes A)  \oplus (A \otimes E )\oplus (E \otimes E)$ 
represents 
$\Delta_A +  \Delta_A^{E,E} + \Delta_E$.

\begin{example} \label{APSexample}
{\rm
Let $A=\Z[X] /(X^2)$ with usual Frobenius algebra structure of truncated polynomial rings,
and $E=\langle Y, Z \rangle_\Z $
with multiplications defined on basis elements by 
$XY=XZ=Y^2 = Z^2=0$ and $YZ=X$.
Comultiplications are defined by
\begin{eqnarray*}
\Delta(1) &=& (1 \otimes X + X \otimes 1)_{A\otimes A}+ (Y \otimes Z + Z \otimes Y)_{E\otimes E}, \\  
\Delta(X) &=& X \otimes X , \\
\Delta(Y) &=& X \otimes Y, \\
\Delta(Z) &=& X \otimes Z.
\end{eqnarray*}
To aid the reader, in the definition of $\Delta(1)$, we indicated
the image of the tensorands. 
The M\"{o}bius  map $\nu$ is defined by 
$\nu (1)= Y + Z$, $ \nu(X)=0$, and $\nu (Y)=\nu (Z)= X$.
This is the structure used in \cite{APS}. 
The correspondence between the above defined structure and their 
simbols are given by  
$1 \leftrightarrow -$,  $X \leftrightarrow +$, 
$Y  \leftrightarrow -0$, $Z  \leftrightarrow +0$.
} \end{example}

\begin{example}\label{TTexample}
{\rm
This is an example that appears 
in \cite{TT}. 
Let $A=E=\Z_2[X,  \lambda^{\pm 1}] /(X^2-\lambda^2 X)$ be the Frobenius algebra
in Example~\ref{univexample} with $h=\lambda^2$,
with $t=0$, and 
 with the coefficient $\Z$ replaced by $\Z_2$.
All mutiplications and comultiplications are those of $A$.
All M\"{o}bius maps are defined by multiplication by $\lambda$.
Then $(A,E)$ gives rise to a commutative Frobenius pair. 

Note that the handle element is $\phi=h=\lambda^2$, so that the M\"{o}bius maps
are multiplication by the square root of the handle element. 
Thus the relation in Fig.~\ref{mobius} (1) and (2) follow immediately, 
and all the others in the figure are automatically satisfied, as the maps are multiplication by a constant.
This definition of the M\"{o}bius maps is derived from  their construction of an unoriented TQFT.
} \end{example}

\begin{remark}\label{ITrem}
{\rm
The above two examples show that commutative Frobenius pairs describe 
the states of categorified Jones polynomials for knots in thickened surfaces 
that are defined 
in \cite{APS,TT}.
We now compare the structures used in \cite{IT} for virtual knots and 
commutative Frobenius pairs. 

Let $A=\Q[X,  t^{\pm 1}] /(X^2-t)$ be the Frobenius algebra
derived from 
Example~\ref{univexample}. 
Thus
\begin{eqnarray*}
\Delta(1) &=& 1 \otimes X + X \otimes 1, \\
\Delta(X) &=& X \otimes X + t \ 1  \otimes 1, 
\end{eqnarray*}
and $\eta(1)=1$. 
Set $E=A$, and define operations as follows.
 Let $\phi $ be the invertible handle element $2X$, so that $ \mu_A \Delta_A=\phi \id_A$. 
 Note that $(2X)^2=4X^2=4t$ 
 is invertible, and so is $\phi$.
 Define  $m^E_{A,E}= \mu_{E,A}^E=   \mu_A$, 
$\Delta_E^{A,E}=\Delta_E^{E,A}=\Delta_A$.
Define further 
$ \mu^A_{E,E} =\phi^{-1}  \mu_A$,  
$\Delta_A^{E,E}=\phi \Delta_A$, 
$\nu_A^E = \phi\  \id$, and $\nu_E^A=\nu_E^E=\id$.
Then $(A,E)$ satisfies most of the axioms of a commutative Frobenius pair, 
but this has the following difference: 
the maps $\mu_E$ and $\Delta_E$ are not defined, but if they were then the identity of Fig.~\ref{samerel} would not be satisfied.

} \end{remark}

\section{TQFTs and commutative Frobenius pairs} \label{TQFTsec}

In this section, we relate commutative Frobenius pairs to 
 topological quantum field theories (TQFTs) of surfaces.  
A $(1+1)$-TQFT is a functor from $2$-dimensional orientable cobordisms to
modules. 
It is 
known that the image of the functor 
 forms a Frobenius algebra~\cite{Abrams,Kock}.
Some aspects  of surface cobordisms  
in $3$-manifolds are studied in \cite{Kaiser}
in relation to TQFTs and surface skein modules~\cite{AF} of $3$-manifolds.
The novel aspects of this paper are to propose 
commutative Frobenius
pairs
for describing TQFTs  in thickened surfaces,  
and to include non-orientable surfaces
by M\"{o}bius maps.
We do not assume that surfaces are orientable throughout this section.
Let $C$ be a properly embedded compact surface
in 
a thickened surface 
$M=F \times I$, where $F$ is a compact 
surface and $I=[0, 1]$. 
We assume that $F$ is not homeomorphic to the projective plane. 
Let $C_0 \sqcup C_1 \subset  C$ be  
the boundary 
$C_i \subset F \times \{ i \}$ for $i=0, 1$,
and $C$ is regarded as a cobordism from $C_0$ to $C_1$. 
Let $\mcob$ be the category of cobordisms of properly embedded 
surfaces  in a  thickened surface $F\times I$
up to ambient isotopy.  
Observe that if $\partial F \ne \emptyset$, then $C_0$ and $C_1$ are embedded in ${\rm{int}} F \times \{0,1\}$. 

We may assume that the height function
$\pi: F \times I \rightarrow I$ restricted to $C$ is a generic Morse
function.
Except at isolated critical levels,
$C_t = \pi^{-1} (t) \cap C$, $t \in [0,1]$, is a set of finite
 simple
 closed curves.
A simple closed curve in $F \times \{ t \}$ is called {\it inessential} 
if it is null-homotopic 
in $F \times \{ t \}$, otherwise {\it essential}. 
Let $C_t=C_t^A \sqcup C_t^E$, $t \in I$,  be the partition of $C_t$ 
into 
inessential curves $C_t^A$ and essential curves $C_t^E$.
In general, for any set $\gamma$ of simple closed curves in $F$,
define $\gamma=\gamma^A\sqcup \gamma^E$ similarly.
Let ${\bf Mod}_k$ be the category of modules over a
commutative unital ring $k$.

\begin{proposition}\label{TQFTthm}
Let ${\cal F}:  \mcob \rightarrow {\bf Mod}_k$ be a TQFT. 
Then the image ${\cal F} (\mcob)$ is a commutative Frobenius pair
with M\"{o}bius maps.  
\end{proposition}
{\it Proof Sketch.\/}
The assignments for surface cobordisms to modules and their homomorphisms 
by a TQFT are 
 made as follows.
Let $C_0=C_t^A \sqcup C_t^E$ be as above for $t=0,1$.
Let $C_0^A =\{ \gamma^A_1, \ldots, \gamma^A_m \}$ and 
$C_0^E =\{ \gamma^E_1, \ldots, \gamma^E_n \}$,
for some non-negative integer $m, n$ (if $m$ or $n$ is $0$, the set is empty). 
Suppose  the functor  
defined on objects $C_0=C_0^A \sqcup C_0^E$
assigns 
$A^{\otimes m}$ to $C_0^A =\{ \gamma^A_1, \ldots, \gamma^A_m\}$
and $E^{\otimes m}$ to $C_0^E =\{ \gamma^E_1, \ldots, \gamma^E_n\}$,
where each 
 tensor factor of $A$ and $E$ is assigned to $\gamma^A_i$ and 
 $\gamma^E_j$ for $i=1, \ldots, m$, $j=1, \ldots, n$,
 hence ${\cal F} (\{ \gamma^A_1, \ldots, \gamma^A_m\} )=A^{\otimes m}$
 and  ${\cal F} (\{ \gamma^E_1, \ldots, \gamma^E_n \} )=E^{\otimes n}$. 
 As in the case of $(1+1)$-TQFTs, saddle points of a cobordism correspond to 
 multiplications and comultiplications, and isotopies correspond to 
 relations such as associativity and coassociativity. 
The rest of the proof is to check that switches of critical levels induce the axioms of a commutative Frobenius pair.

\begin{figure}[htb]
\begin{center}
\includegraphics[width=3.5in]{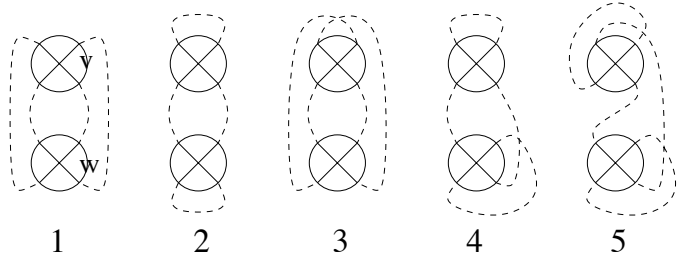}
\end{center}
\caption{The list of two-crossing connections from Reference \cite{APS}}
\label{APSlist}
\end{figure}

 The axioms are checked using cobordisms with generic Morse functions,
using the list in \cite{APS}. We illustrate 
the gist of the 
argument. 
In Fig.~\ref{APSlist}, we copied the list in \cite{APS} of all possible connections
of a pair of crossings for link diagrams in $F \times I$. 

\begin{figure}[htb]
\begin{center}
\includegraphics[width=4in]{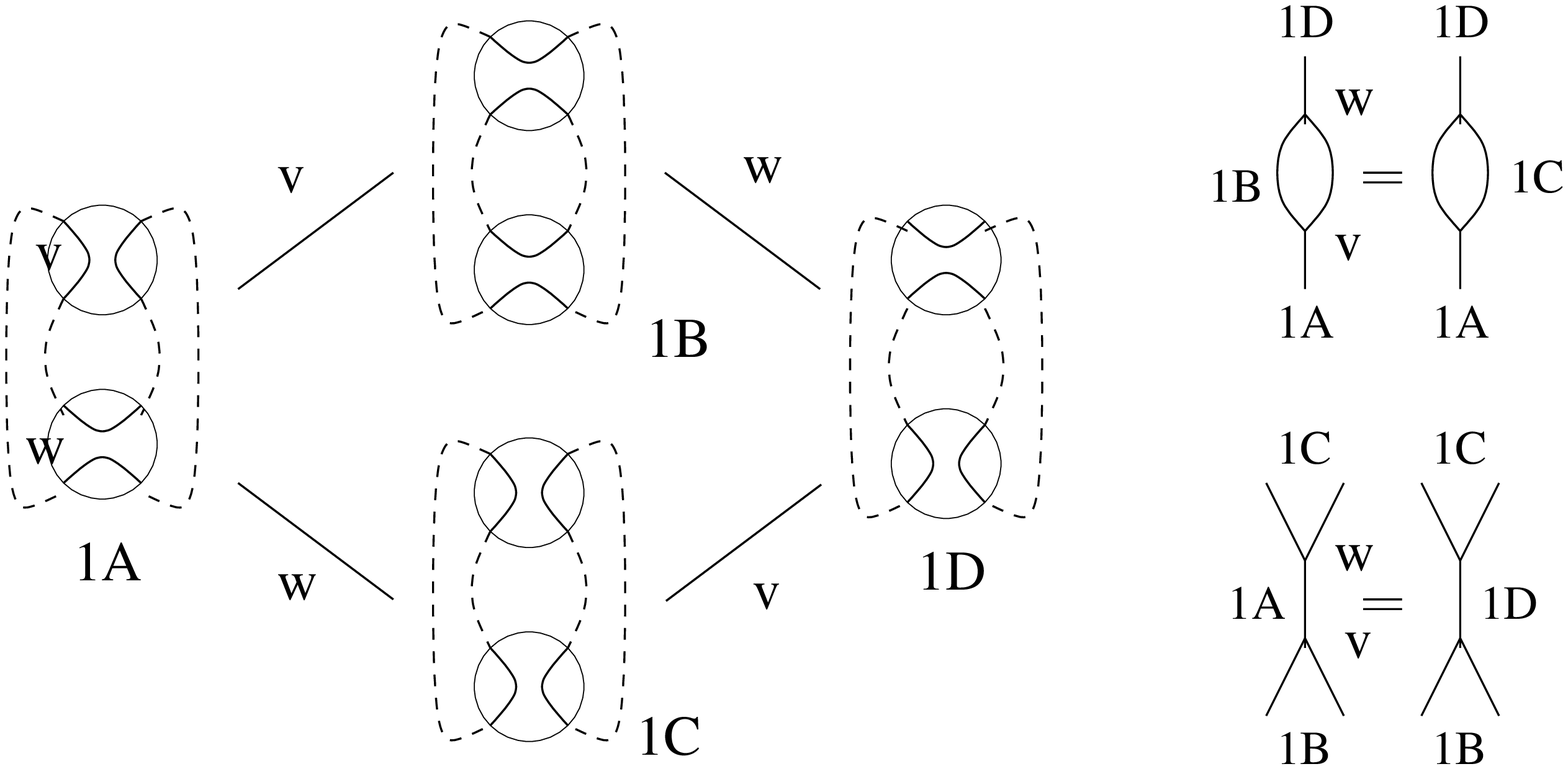}
\end{center}
\caption{Case 1 from Reference \cite{APS}}
\label{case1}
\end{figure}

The crossing can be viewed as a non-generic cross section at a saddle point. Thus the diagram 1 of Fig.~\ref{APSlist} can be thought of as being expanded to the diamond illustrated in Fig.~\ref{case1}. There are four possible generic perturbations of this diagram that are indicated by the vertices of the diamond (1A through 1D). The edges of the diamond represent passing through the corresponding saddle point. 
The positively sloped edges are passing through the upper saddle point $v$, and the negatively sloped edges are passing through the lower saddle point 
$w$.

In the right of Fig.~\ref{case1}, diagrams  are depicted 
that represent
 the maps corresponding to the cobordisms of surfaces.  For example, the top left diagram starts with $1A$ that has a single foot, 
then separates into two arcs and merges to a single arc, 
representing $1A$ is connected, $1B$ has two components, and $1D$ has one component.
Then as a cobordism this is equal to the cobordism corresponding to $1A-1C-1D$, 
which is depicted in the right of the equality.

If the curve in Figure $1A$ is inessential, then the two curves  after the first smoothing along the
 $(ABD)$-path are either both inessential or both essential. So the $(ABD)$-path  maps to either the composition   $\mu_A \circ \Delta_A$ or   the composition $\mu^A_{E,E} \circ \Delta^{E,E}_A.$ Meanwhile, $(ACD)$-path is the opposite composition. In this case, the equality in Fig.~\ref{mobius} (2) applies. 

If the curve in Figure $1A$ is essential, it either splits into two essential curves and these merge to an essential curve, 
 or it splits into an essential curve and an inessential curve that merge to an essential curve. In this case, the identity of Fig.~\ref{samerel} (3) applies. 

Similarly, we compare the paths $(BAC)$ and $(BDC)$, 
 and test the equality represented in the bottom right of Fig.~\ref{case1}. There are many cases to consider: both feet of the $1B$ state are essential $(E \otimes E),$ both are inessential $(A\otimes A)$ or a mixed state $(E \otimes A)$. After the merger, the center curve can be either essential or inessential. Most often the equality holds by default, however the case represented in Fig.~\ref{samerel} (2) can also occur. 
 
The condition that $F$ is not a projective plane
is necessary from the following fact: 
when a connected circle goes through a saddle and becomes another connected 
 curve, at least one of them (before or after the saddle) must be essential
 under the assumption that 
$F$ is not  a projective plane.

The rest of the proof involves considering all the possible cases represented in Fig.~\ref{APSlist}, comparing paths in the corresponding smoothing square, and examining the cases among essential and inessential curves. 
All the  conditions are satisfied by axioms, 
or their consequences such as those in Fig.~\ref{maplemmas}.
$\Box$

\begin{figure}[htb]
\begin{center}
\includegraphics[width=2in]{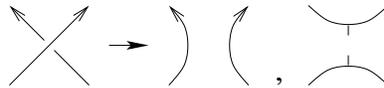}
\end{center}
\caption{Smoothings for  the Miyazawa polynomial}
\label{Miyazawa}
\end{figure}

\begin{figure}[htb]
\begin{center}
\includegraphics[width=5in]{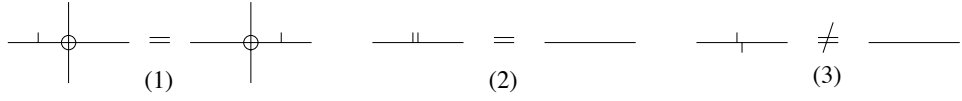}
\end{center}
\caption{Poles in the Miyazawa polynomial}
\label{poles}
\end{figure}

In the remaining of this section, we relate commutative Frobenius pairs to
surface cobordisms with poles that appear in Miyazawa polynomials~\cite{Miya}
and extended bracket~\cite{Lou} for virtual knots. 
In the skein relation of the Miyazawa polynomial \cite{Miya}
(see also  \cite{Lou})
  smoothings of crossings with {\it poles} 
  were used. 
  When
 the smoothing that does not respect orientations of arcs 
 is performed, 
 a pair of short segments, called poles,   
 are placed at arcs after the smoothing
 (see Fig.~\ref{Miyazawa}).
 In defining  their invariants, the relations satisfied between poles and virtual crossings 
 are depicted in Fig.~\ref{poles} (1) and (2). In (3), a relation that is not imposed is depicted.
 
\begin{definition} {\rm 
A set of  {\it virtual circles with poles} 
is a virtual link diagram on the plane without classical (over-under) crossings,
with an even number (possibly zero) of poles attached.

Two sets of virtual circles with poles are {\it equivalent } if
they are related by a finite sequence of virtual 
Reidemeister moves without classical crossings, 
cancelation/creation of a pair of adjacent poles on the same side
as 
depicted in Fig.~\ref{poles} (2), and plane ambient isotopies. 
} \end{definition}

\begin{definition} {\rm 
A pair of adjacent poles on opposite sides are not canceled as depicted in Fig.~\ref{poles} (3),
and the number of such pairs on a transverse component $C$, 
after canceling poles on the same side as depicted in 
Fig.~\ref{poles} (2),   
is called 
the {\it degree}, and denoted by $\vd{ (C) }$.

The sum of the degrees of all components of a set ${\cal V} $ of virtual circles with poles
is called the degree of ${\cal V} $ and denoted by $\vd({\cal V} )$.
The degree is well-defined up to equivalence of 
virtual circles with poles.

A set of virtual circles with poles is 
{\it essential} if $\vd({\cal V} )>0$, inessential otherwise. 
} \end{definition}

\begin{figure}[htb]
\begin{center}
\includegraphics[width=4in]{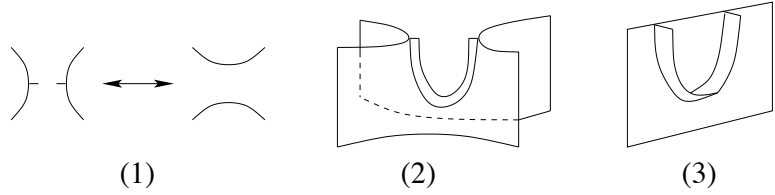}
\end{center}
\caption{A saddle guided by poles}
\label{polesaddle}
\end{figure}

\begin{definition} {\rm
Two  sets of virtual circles with poles are {\it related by a poled saddle}
if they are identical except a small disk neighborhood in which 
they are as depicted in the left and right of Fig.~\ref{polesaddle} (1). 

Two  sets of virtual circles with poles are {\it pole cobordant}
if they are related by a {\it pole cobordism}: a finite sequence of 
equivalences, poled saddles, and birth/death of trivial inessential curves.
} \end{definition}

It is known~\cite{CS:book} that Reidemeister moves 
are derived from  cross sections of surface cobordisms of generic surfaces.
Thus a pole cobordism is regarded as a surface cobordism of 
a generic surface in $\R^2 \times [0,1]$ with 
the 
continuous 
images of poles 
added along saddles, {\it etc.}. 
These vestigial surfaces are called {\it hems}.
This is expressed in Fig.~\ref{polesaddle} (2) for a saddle point corresponding to
a poled saddle. There is a  poled saddle with the upside-down picture. 
The hem is always on the negatively curved side of  a saddle. 
In Fig.~\ref{polesaddle} (3), 
a cancelation of a pair of poles on the same side 
that 
corresponds to a minimal point of hem 
is depicted. 

\begin{definition} {\rm
Let the category of virtual cobordisms $\vcob$ be the category with 
the objects the finite sets of   virtual circles with poles, and the
morphisms 
generated by 
the finite sets of  pole cobordisms.
} \end{definition}

\begin{proposition}\label{virtTQFTthm}
Let  ${\cal F}' : \vcob \rightarrow {\bf Mod}_k$ be a TQFT.
Then the image ${\cal F}'  ( \vcob ) $ is a commutative Frobenius pair
with M\"{o}bius maps.
\end{proposition}
{\it  Sketch Proof.\/} 
Assignments  are made from virtual circles with poles to
modules $A$ and $E$   
in the same way as in Proposition~\ref{TQFTthm}.  
Then axioms are checked for all cases of \cite{APS} as in the proof of Proposition~\ref{TQFTthm}. 
$\Box$

\bigskip

Let $(A, E)$ be a commutative Frobenius pair over $k$.
Let {\bf FP}$(A, E)$ be the subcategory of ${\rm Mod}_k$
whose objects are generated by $A^{\otimes m} \otimes E^{\otimes n}$
for non-negative integers $m$ and $n$, and
morphisms generated by the maps
$m, \Delta, \eta, \epsilon, \tau$ of the commutative Frobenius pair $(A, E)$,
and its M\"{o}bius maps. 

We conjecture that 
 there are functors 
${\cal F}: \mcob \rightarrow \mbox{\bf FP} (A,E)$ and 
${\cal F}' : \vcob  \rightarrow \mbox{\bf FP} (A,E)$ as defined  in the proofs of 
Propositions~\ref{TQFTthm} and \ref{virtTQFTthm},
respectively. 

Typically,  
to prove such a conjecture it is shown 
that the stated set of relations is sufficient 
to describe isotopy  
which, in turn, is proved by using relations to deform a given cobordism 
to a standard form.

\section{Constructions of commutative Frobenius pairs} \label{constructsec}

We consider commutative Frobenius pairs $(A, E)$
 for 
$A=\Z[X, h, t]/(X^2-hX-t)$,
as this appears in the universal Khovanov homology \cite{Kh06}.
See  Example~\ref{univexample} for the Frobenius algebra structure.

Suppose $E$ is of rank $2$, and denote by $Y,Z$ the basis elements of $E$.
Motivated from Example~\ref{APSexample}, 
we characterize the case when the multiplication on $E$ is trivial:
  $\mu^E_{E,E}=\Delta^{E,E}_E=0$.

\begin{theorem} \label{APSgenThm}
Let $A=\Z[X, h, t ]/(X^2-hX-t)$ and $E=\langle Y, Z \rangle$.
If $(A, E)$ is a commutative Frobenius pair
with $\mu^E_{E,E}=\Delta^{E,E}_E=0$, then $A$ must be  of the form 
$A=\Z[X, a ]/ (X-a)^2 $,
so that $h=2a$ and $t=-a^2$. 

Let  $A=\Z[X, a ]/(X-a)^2 $,  
$E=\langle Y, Z \rangle$ and assume  $\mu^E_{E,E}=\Delta^{E,E}_E=0$.
Then 
 $(A,E)$ is a commutative Frobenius pair with M\"{o}bius maps 
 if and only if the following conditions are satisfied.
 The multiplications and actions are defined by 
$$ XY=aY, \ XZ=aZ, \ Y^2=c_{YY} (X-a), \ YZ= c_{YZ} (X-a), \ Z^2=c_{ZZ} (X-a), $$ 
for some constants $c_{YY}$, $c_{YZ}$, and $c_{ZZ}$ that satisfy the conditions below. 
 The comultiplications and coactions are defined by 
 \begin{eqnarray*}
 \Delta_E^{A,E} (Y)&=& (X-a) \otimes Y, \\ 
\Delta_E^{A,E} (Z)&=&  (X-a) \otimes Z, \\  
\Delta_A^{E,E}(1) &=& d_{YY} Y \otimes Y + d_{YZ} Y \otimes Z 
+d_{YZ} Z \otimes Y +  d_{ZZ} Z \otimes Z , \\
\Delta_A^{E,E}(X) &=& a\  \Delta_A^{E,E}(1) ,  
\end{eqnarray*}
for some constants $d_{YY}$, $d_{YZ}$, $d_{ZZ}$ that satisfy the conditions below.

The M\"{o}bius maps are defined by 
$$\nu_E^E=0, \quad \nu_E^A(Y)=e_{Y}(X-t), \quad \nu_E^A(Z)=e_{Z}(X-a),  
\quad \nu_A^E(1)=f_Y Y + f_Z Z . $$
The coefficients satisfy the following conditions.
Let  $C$ be the $2 \times 2$ symmetric  integral  matrices with entries 
 $c_{ij} $ with $\{ i, j \} =\{ Y, Z \} $, respectively, and let 
 $\vec{e}=[ e_{Y} ,  e_{Z} ]^T$,  
$\vec{f}=[  f_{Y} ,  f_{Z}]^T$ be column vectors.
Then they satisfy
 $$C \vec{f}=\vec{e},  \quad \vec{e} \cdot \vec{f}=2, \quad 
c_{YY} d_{YY}+ 2  c_{YZ} d_{YZ} + c_{ZZ}  d_{ZZ} =2 \quad \pmod{(X-a)}. $$  
\end{theorem}
{\it Proof.} 
First we prove the following lemmas.

\begin{lemma}
\label{first}
 Suppose that $A=\Z[X, h, t]/(X^2-hX-t)$, and that $E=\langle Y,Z \rangle$ is a left module for $A$,
 where $k$ is an integral domain. 
  Let 
 $\mu^E_{A,E} (X \otimes Y) = XY=a_0 Y + a_1 Z$, and 
$\mu^E_{A,E} (X \otimes Z) = XZ=b_0 Y + b_1 Z$.
Then 
\begin{eqnarray*}
 a_0^2 + a_1 b_0  -a_0 h - t =0, & & 
  (a_0  +  b_1 -h ) a_1 =0,  \\ 
a_1 b_0 + b_1^2 - h b_1 - t =0 ,  & & 
(a_0  +  b_1 - h )b_0  = 0.
\end{eqnarray*}
\end{lemma}
{\it Proof.} The first two equations follow from the equation $X(XY)= X^2 Y$, and the last two follow from $X(XZ)=(X^2)Z.$ $\Box$

\begin{remark}{\rm 
For 
a commutative Frobenius pair $(A, E)$ 
as in Lemma~\ref{first},    
the conditions 
$\Delta^{A,E}_E = (\id_A \otimes \mu^E_{A,E}) (\Delta_A(1)\otimes \id_E) $ and
$\Delta_E^{A, E}  \mu_{A,E}^E=( \mu_{A,A}^A \otimes  \id_E )(\id_A \otimes \Delta_E^{A,E})$
applied to  $X \otimes Y$ and 
$X \otimes Z$  to give the same relations.} \end{remark}

\begin{lemma}\label{htzero}
\begin{sloppypar}
Suppose that the commutative Frobenius pair $(A,E)$ satisfies the hypotheses of Lemma~\ref{first}, and suppose further that 
$\mu^E_{E,E}=\Delta^{E,E}_E=0$.
Then 
there are constants $a$, $c_{YY}$, $c_{YZ}$ and  $c_{ZZ}$ such that 
$h=2a$, $t=-a^2$ and the following hold:
$XY=aY$, $XZ=aZ$, 
$Y^2=c_{YY} (X-a)$, $YZ=c_{YZ} (X-a)$ and $Z^2=c_{ZZ} (X-a)$. 
\end{sloppypar}
\end{lemma}
{\it Proof.\/} 
First we show that  $h^2 + 4t =0$, and $XY=aY$, $XZ=aZ$, where 
$h=2a$, $t=- a^2$. 
The assumption and the relation depicted in Fig.~\ref{samerel} (3) implies 
$ \mu_{A,E}^E \Delta_E^{A, E} = \mu_{E,E}^E \Delta_E^{E, E}=0$. 
On the other hand, 
$$ \mu_{A,E}^E \Delta_E^{A, E} (Y)
=  \mu_{A,E}^E(\id_A \otimes \mu^E_{A,E}) (\Delta_A(1)\otimes \id_E)(Y)=2XY-hY ,$$
so that we have $2(a_0 Y + a_1 Z) - hY=0$, and similarly for $Z$ we obtain
$2(b_0 Y + b_1 Z) - hY=0$. 
Thus we have $2a_0=h=2b_1$, $2a_1=0 = 2 b_0$. 
It follows that $a_1=b_0=0$, and from $a_0^2 + a_1 b_0  -a_0 h - t =0$ 
in Lemma~\ref{first}, we have 
 $a_0^2 +t=0=b_1^2 +t$.
 Set $a=a_0=b_1$, and we obtain the result.

Set  $XY=aY$ and $XZ=aZ$, where $h=2a$ and $t=- a^2$ as above.
Set  $Y^2=b_{YY} \ 1_A + c_{YY} X$. 
Here we abused the notation for
$Y^2$ to be $ \mu_{E,E}^A(Y\otimes Y) $ as we assumed   $\mu^E_{E,E}=0$.
We compute 
\begin{eqnarray*}
(XY)Y &=& aY^2= a(b_{YY} \ 1_A + c_{YY} X) \\
X(Y^2)&=& X(b_{YY} \ 1_A + c_{YY} X)= b_{YY} X  + c_{YY} (hX + t) \\
 & = & (b_{YY}  + c_{YY} h)X + c_{YY} t, 
 \end{eqnarray*}
 hence $b_{YY}  + c_{YY} (h-a)=b_{YY}  + c_{YY} a=0$ and 
  $a b_{YY}-  c_{YY} t =a(b_{YY} + c_{YY}a) =0$. 
This implies 
$ a=0$ or $b_{YY}=-a c_{YY}$. 
Set $YZ=b_{YZ} \ 1_A + c_{YZ} X$ and $Z^2=b_{ZZ} \ 1_A + c_{ZZ} X$.
A similar argument for $XZ^2$ and $XYZ$ shows that 
[ $a=0$  or  $b_{YZ}=-a c_{YZ}$ ], and [ $a=0$  or  $b_{ZZ}=-a c_{ZZ}$ ], respectively.
If $a=0$, then $h=t=0$ and  $b_{YY}=b_{YZ}=b_{ZZ}=0$.
If $a\neq 0$, then $Y^2=c_{YY} (X-a)$, $YZ=c_{YZ} (X-a)$ and $Z^2=c_{ZZ}(X-a)$. 
Either way the result follows.
$\Box$
 
\bigskip 

\noindent
{\it Proof (of Theorem~\ref{APSgenThm}) continued.\/} 
First we determine comultiplications.
{}From Lemma~\ref{htzero}, we have 
$h=2a$, $t=-a^2$, 
$XY=aY$, $XZ=aZ$, 
$Y^2=c_{YY} (X-a)$, $YZ=c_{YZ} (X-a)$ and $Z^2=c_{ZZ} (X-a)$. 
Note that $X^2-hX -t=(X-a)^2$. 
Since $\Delta_A^{E,E} m_{E,E}^A(Y \otimes Y)=\Delta_E^{E,E} m_{E,E}^E(Y \otimes Y)=0$
(see Fig.~\ref{samerel} (2)), 
we obtain $c_{YY} \Delta_A^{E,E}(X-a)=0$. 
If $ \Delta_A^{E,E}(X-a) \neq 0$, then 
$Y^2=YZ=Z^2=0$, and the equality 
$$  \mu_{E,E}^A \Delta_A^{E,E} (1) =  \mu_{A,A}^A \Delta_A^{A,A} (1)=2X-h=2(X-a)$$
 in Fig.~\ref{mobius} (2) leads to a contradiction. 
 Hence $ \Delta_A^{E,E}(X-a)=0$, so that $ \Delta_A^{E,E}(X)=a \Delta_A^{E,E}(1)$. 
 Note also that at least one of $c_{YY}$, $c_{YZ}$ and $c_{ZZ}$ is non-zero. 
  {}From  $\Delta_E^{A, E} (Y)
= (\id_A \otimes \mu^E_{A,E}) (\Delta_A(1)\otimes \id_E)(Y) $
(see Fig.~\ref{cancel}), 
we obtain $\Delta_E^{A,E}(Y)=(X-a) \otimes Y$, and similarly,
$\Delta_E^{A,E}(Z)=(X-a) \otimes Z$. 
Let 
$$\Delta_A^{E,E} (1)= d_{YY} Y \otimes Y + d_{YZ} Y \otimes Z +d_{YZ} Z \otimes Y +d_{ZZ} Z \otimes Z.$$
Here we used the cocommutativity.
Then from 
Fig.~\ref{mobius} (2), 
we obtain 
$$ (c_{YY} d_{YY}+ 2  c_{YZ} d_{YZ} + c_{ZZ}  d_{ZZ} )  (X-a)  = 2(X-a) . $$
Hence $c_{YY} d_{YY}+ 2  c_{YZ} d_{YZ} + c_{ZZ}  d_{ZZ} =2+c (X-a) $ for some constant $c$,
i.e.,  $c_{YY} d_{YY}+ 2  c_{YZ} d_{YZ} + c_{ZZ}  d_{ZZ} =2$ modulo $(X-a)$.

Next we determine M\"{o}bius maps. 
{}From  
Fig.~\ref{mobius} (7), 
we have $\nu_E^E=0$. 
{}From the up-side down diagram of 
Fig.~\ref{mobius} (4)  
with input $1 \otimes X$, we have 
$\nu_A^E (X)=X \nu_A^E (1)  = a \nu_A^E (1)$ since $XY=aY$ and $XZ=aZ$. 
Set $\nu(1)=f_Y Y + f_Z Z$.
{}From  the up-side down diagram of 
Fig.~\ref{mobius} (5)  
with input $X \otimes Y$, 
we have  
$ X \nu_E^A(Y) = \nu_E^A (XY) =  \nu_E^A ( aY) = a \nu_E^A (Y)$,
so that $ \nu_E^A(Y)  (X-a)=0$, and we obtain $ \nu_E^A(Y) =e_Y (X-a)$ 
for some constant $e_Y$. Similarly, 
 $ \nu_E^A(Z) =e_Z (X-a)$ 
{}From Fig.~\ref{mobius} (1) with input $1$, we have 
 $(e_Y f_Y + e_Z f_Z) (X-a) =2 (X-a)$,  
hence  $\vec{e} \cdot \vec{f} = 2 $ modulo $(X-a)$.  
The up-side down diagram of Fig.~\ref{mobius} (6) with inputs 
$1 \otimes Y$ and $1 \otimes Z$ gives 
$C \vec{f}=\vec{e}$ modulo $(X-a)$. 
 
Thus the conditions stated in this theorem are necessary, and 
it is checked that they are also sufficient.
$\Box$

\begin{remark}
{\rm 
Example~\ref{APSexample} is the case when 
$\vec{e}=\vec{f}=[1,1]^T$, 
$c_{YY}=c_{ZZ}=d_{YY}=d_{ZZ}=0$, $c_{YZ}=d_{YZ}=1$.
}
\end{remark}

\begin{theorem} \label{rootThm}
Let $A$ be a commutative Frobenius algebra 
over a commutative unital ring $k$ 
with handle element $\phi$, 
such that there exists an element $\xi \in A$ with $\xi^2=\phi$.
Then there exists a commutative Frobenius pair with M\"{o}bius maps.
\end{theorem}
{\it Proof.\/}
Define all multiplications, comultiplications, action and coaction  by those of $A$.
Then all axioms for these operations 
are satisfied, and only the M\"{o}bius maps remain to be defined and 
their axioms checked.
Define all the M\"{o}bius maps by multiplication by $\xi$.
The 
relations (1) and (2) in Fig.~\ref{mobius}
involving M\"{o}bius maps and the handle element $\phi$
are satisfied since $\xi^2=\phi$.
Other relations in Fig.~\ref{mobius}
 involving M\"{o}bius maps and (co)multiplications are satisfied since
the former are assigned multiplication by a constant.
This construction was inspired by Example~\ref{TTexample}.
$\Box$

\begin{example}{\rm
One of the examples of Frobenius algebras used in \cite{TT}
is  $k=\F_2 [\lambda] $ and $A=k[X]/(X^2 - \lambda^2 X)$.
The handle element is indeed $\lambda^2$, and we can take $\xi=\lambda$. 
}
\end{example}

\begin{corollary} \label{rootcor}
Let $A=E=k[X]/(X^2-hX-t)$, where $k=\Z[a^{\pm 1}, b^{\pm1}]$, 
and  $h=-2b^{-1} (a - b^{-1})$, $t=-b^{-2}(a^2 + h).$
Then  $(A,E)$ gives rise to  a commutative Frobenius pair with M\"{o}bius maps.
\end{corollary}
{\it Proof.\/}
Let $\xi=a + bX$, and one computes that $\xi^2=\phi=2X-h$.  
$\Box$

\begin{figure}[htb]
\begin{center}
\includegraphics[width=4.5in]{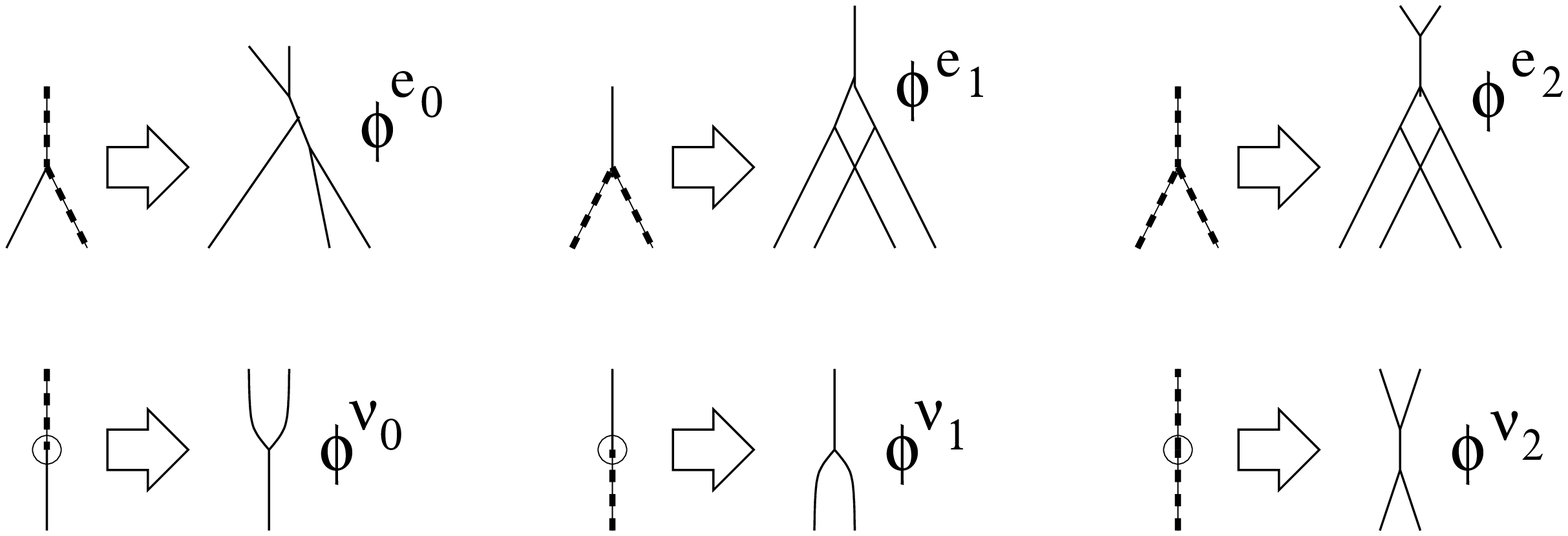}
\end{center}
\caption{Double tensor assignments on $E$}
\label{dbles}
\end{figure}

\begin{theorem} \label{dbleThm}
If $A$ is a  commutative Frobenius algebra 
over a commutative unital 
ring  $k$, such that its handle element $\phi \in A$
is invertible,  then 
 there exists a commutative Frobenius pair $(A,E)$ with M\"{o}bius maps.
\end{theorem} 
{\it Proof.\/} Let $E=A \otimes A$.
Then define various multiplications, actions and M\"{o}bius maps as depicted in 
Fig.~\ref{dbles}. In the figure, powers of $\phi$ are indicated for  assignments, 
whose exponents are specified below. Each assigned map is the map indicated 
by the diagrams multiplied by the specified power of $\phi$. 
Specifically,
they are defined by
\begin{eqnarray*}
\mu_{A,E}^E &=& \phi^{e_0} \Delta_A  \mu_{A} (\id \otimes  \mu_A ), \\
 \mu_{E,E}^A &=&  \phi^{e_1}  \mu_A ( \mu_A \otimes  \mu_A) (\id \otimes \tau \otimes  \id) , \\
 \mu_{E,E}^E &=&  \phi^{e_2}  \Delta_A \mu_A ( \mu_A \otimes  \mu_A) (\id \otimes \tau \otimes  \id) , \\
 \nu_A^E &=& \phi^{\nu_0} \Delta_A, \\
  \nu_E^A &=& \phi^{\nu_1} \mu_A, \\
   \nu_E^E &=& \phi^{\nu_2} \Delta_A \mu_A .
\end{eqnarray*}

The assignment for comultiplications and coactions 
are defined by up-side down diagrams of multiplications without 
 $\phi$ factors. 
Then the axioms of the commutative Frobenius pairs are 
satisfied with the unique choice $e_0=-1$, $e_1=e_2=-2$,
$\nu_0=1$, $\nu_1=-1$ and $\nu_2=0$. 
The assignment of $\phi$ factors is inspired by construction in \cite{IT} as described in 
Remark~\ref{ITrem}.
$\Box$

\begin{corollary}
For $A=k[X]/(X^2 -hX - t ) $, where $h, t$ are variables taking values in $k$, 
if $4t + h^2 \in k$ is invertible, then 
 there exists a commutative Frobenius pair $(A,E)$ with M\"{o}bius maps.
 \end{corollary}
{\it Proof.\/}
The handle element is computed as
$\phi = m \Delta (1)=2X -h$. 
Then one computes $\phi^2=4t+ h^2$, which is assumed to be invertible, 
so that $\phi ^{-1} = (4t + h^2)^{-1} \phi$, 
and the result follows from 
Theorem~\ref{dbleThm}.
$\Box$

\subsection*{
Acknowledgments}

JSC was supported in part by NSF Grant DMS \#0603926. MS was supported in
part by NSF Grants DMS \#0603876.

\end{document}